%% file: arxiv.tex
\RequirePackage{fix-cm}
\documentclass[smallcondensed]{svjour3}     
\smartqed  
\usepackage{graphicx}
\usepackage{mathptmx}      
\usepackage{inputenc, theorem,amsmath,enumerate,fancyhdr,amssymb,amsfonts,multirow, url, float}
\usepackage[round]{natbib}

\newcommand{\norm}[1]{\left\|{#1}\right\|}
\newcommand{\pp}{{p}}
\newcommand{\bigO}{\mathcal{O}}
\newcommand{\conv}{\operatorname{conv}}
\newcommand{\ri}{\operatorname{ri}}
\newcommand{\aff}{\operatorname{aff}}

\newcommand{\Sub}{\operatorname{Sub}}

\newcommand{\Null}{\operatorname{Null}}
\newcommand{\Span}{\operatorname{Span}}

\newcommand{\RR}{\mathbb{R}}

\newcommand{\C}{\mathcal{C}}

\newcommand{\J}{\mathcal{J}}

\newcommand{\Pp}{\mathcal{P}}

\newcommand{\Ss}{\mathcal{S}}

\begin{document}

\title{A faster dual algorithm for the Euclidean minimum covering ball problem
}

\titlerunning{}        

\author{Marta Cavaleiro        \and
        Farid Alizadeh
}


\institute{Marta Cavaleiro\at
              Rutgers University, MSIS Dept. \& RUTCOR \\
              100 Rockafeller Rd, Piscataway, NJ 08854.
              \email{marta.cavaleiro@rutgers.edu}           
           \and
           Farid Alizadeh \at
              Rutgers University, MSIS Dept. \& RUTCOR \\
              100 Rockafeller Rd, Piscataway, NJ 08854.
              \email{farid.alizadeh@rutgers.edu} 
}

\date{Received: date / Accepted: date}

\maketitle

\begin{abstract}
Dearing and Zeck \citeyearpar{Dearing09} presented a dual algorithm for the problem of the minimum covering ball in $\mathbb{R}^n$. Each iteration of their algorithm has a computational complexity of at least $\bigO(n^3)$. In this paper we propose a modification to their algorithm that, together with an implementation that uses updates to the QR factorization of a suitable matrix, achieves a $\bigO(n^2)$ iteration. 

\end{abstract}

\section{Introduction}
\input{1Introduction}

\section{The dual algorithm by Dearing and Zeck}\label{sec:DearingZeck}
\input{2DearingZeck}

\section{Directional search via orthogonal projections}\label{sec:LineSearch}
\input{3NewLineSearch}

\section{The Implementation}\label{sec:Implementation}
\input{4Implementation}

\section{Computational results}\label{sec:Results}
\input{5Results}



\bibliographystyle{spbasic}      
\bibliography{referencesDualAlg}

\end{document}

%% file: 1Introduction.tex
Consider a given set of points $\Pp=\{p_1,\hdots,p_m\}$ in the Euclidean space $\RR^n$. Let~$\|.\|$ denote the Euclidean norm. The problem of finding the hypersphere $$B(x, r)=\{y\in\RR^n: \|y-x\|\leq r\}$$ with minimum radius that covers $\Pp$, which we will refer to as the minimum covering ball (MB) of $\Pp$, can be formulated as

\begin{equation} MB(\Pp) := \begin{array}{cl}\label{meb}
\min &r^2 \\
\text{s.t.} & \|p_i-x\|^2\leq r^2, \quad i=1,\hdots,m.
\end{array}
\end{equation}

\noindent We will use the notation $MB(\Pp)$ both to refer to the problem of the minimum covering ball of a set $\Pp$ and, depending on the context, also to the corresponding optimal ball.

\medskip

The MB problem, reported to date back to the 19th century \citep{Sylvester}, is an important and active problem in computational geometry and optimization. Applications include facility location, see e.g. \cite{Hale03}, \cite{Moradi09}, and \cite{Plastria02}; computer graphics, see e.g. \cite{Hubbard96}, and \cite{Larsson16}; machine learning, see e.g. \cite{Kumar03}, \cite{Nielsen09}, and the references therein; etc.

Problem (\ref{meb}) can easily be converted in a quadratic program (QP) and solved using off-the-shelf QP solvers. G{\"a}rtner and Sch{\"o}nherr \citeyearpar{Gartner00} developed a generalization of the simplex method for QP with the goal of targeting geometric QPs, with one of the main applications being the MB problem. The MB problem is also a second-order-cone program (SOCP) instance, so interior point methods may also be used. Some work has also been done using the SOCP formulation while exploiting the special features of the MB problem \citep{Kumar03, Zhou05}. 

The algorithmic complexity of the MB problem was first given by Megiddo in 1983, when he gave the first deterministic algorithm that solves the MB problem in linear time on the number of points, when the dimension is fixed \citep{Megiddo83,Megiddo84}.

The MB problem is an \emph{LP-type} problem \citep{Dyer04}, sharing common properties with linear programs, in particular its combinatorial nature: the optimal ball is determined by at most $n+1$ points of $\Pp$ that lie on its boundary. Such set of points is usually called a \emph{support set}. Many algorithms that search for such set have been developed. Welzl \citeyearpar{Welzl91} presented a randomized algorithm that searches for a support set, solving the problem in expected linear time for fixed dimension. This algorithm was improved by G\"artner \citeyearpar{Gartner99}, however only dimensions $n\leq 25$ could be handled in practice. Later Fischer and G{\"a}rtner \citeyearpar{Fischer04} proposed an algorithm with a pivoting scheme resembling the simplex method for LP based on previous ideas from \citep{Hopp96}, even adapting Bland's rule to avoid cycling. Their algorithm can handle problems in larger dimensions ($n\leq 10,000$). Using related ideas, Dearing and Zeck \citeyearpar{Dearing09} developed a dual algorithm for the MB problem. This algorithm is the subject of this paper, and as we will prove, it can also deal with larger dimensions.



The concept of $\epsilon$-\emph{core sets}
, proposed by B\^{a}doiu et al. \citeyearpar{Badoiu02}, introduced a new direction of research in approximate algorithms for the MB problem. A remarkable property is the existence of a $\epsilon$-core set of size at most $1/\epsilon$ and independent of $n$ \citep{Badoiu03,Kumar03}. Several algorithms focused on finding $\epsilon$-core sets have been proposed \citep{Badoiu02, Badoiu03, Kumar03, Yildirim08, Nielsen09, Larsson13}, being in general able to deal with large dimensions in useful time.

Streaming algorithms, that only allow one pass over the input points, have also been studied. Zarrabi-Zadeh and Chan \citeyearpar{Zarrabi06} gave a $3/2$-approximation algorithm, and later an algorithm by Agarwal and Sharathkumar \citeyearpar{Agarwal15} was able to achieve a $1.22$-approximation factor \citep{Chan11}.

\medskip

In this paper we propose a modification to the algorithm proposed by Dearing and Zeck \citeyearpar{Dearing09} that makes it faster. Their algorithm looks for a support set, by solving a sequence of subproblems $MB(\Ss)$, with $\Ss\subseteq\Pp$ affinely independent, until all points of $\Pp$ are covered. At each iteration, set $\Ss$ is updated by either adding a point that is not yet covered, or by replacing an existing point in $\Ss$ by it. Problem $MB(\Ss)$ is solved using a directional search procedure, and during this step possibly more points are removed from $\Ss$. 

It is possible to implement the algorithm as presented in \citep{Dearing09} taking advantage of the QR updates of a suitable matrix, which can be done in quadratic time on $n$. However, they would need to be done as many as $\bigO(n)$ times, resulting in an iteration having $\bigO(n^3)$ computational complexity. We modify the directional search procedure in such a way that, together with an implementation using QR updates, achieves a $\bigO(n^2)$ iteration.

\smallskip

The paper is organized as follows: In section~\ref{sec:DearingZeck} we review Dearing and Zeck's algorithm for the $MB$ problem. In section~\ref{sec:LineSearch} we present a modification of the algorithm that makes its directional-search step faster. In section~\ref{sec:Implementation} we show how the algorithm can be implemented using QR updates of a certain matrix, which together with the results from section~\ref{sec:LineSearch}, result in $\bigO(n^2)$ iteration. Finally, in section~\ref{sec:Results} we present some computational results that show the practical impact of our work.

\smallskip

One application of this algorithm that benefits from it having a faster iteration is a branch and bound approach to solve the problem of the minimum $k$-covering ball (that seeks the ball with smallest radius that contains at least $k$ of $m$ given points in $\Pp$). At each node of the search tree, Dearing and Zeck's algorithm can be employed to solve the corresponding subproblem. Since it starts with the solution of the parent node, it usually needs very few iterations to solve each subproblem. And since a very large number of subproblems need to be solved, having a fast iteration is essential. This application will be the subject of a future publication by the authors of this paper.

%% file: 2DearingZeck.tex
The dual problem of (\ref{meb}) is
\begin{equation} \begin{array}{cl}\label{dual}
\displaystyle\max_{\pi\in \RR^m} & \displaystyle\sum_{i=1}^{m} \pi_i\left\|p_i\right\|^2 - \left\|\sum_{i=1}^{m} \pi_i p_i\right\| ^2\\
\text{s.t.} & \displaystyle\sum_{i=1}^{m} \pi_i =1\\
& \pi_i\geq 0,\quad i=1,\dots, m,
\end{array}\end{equation}
\noindent and the optimal primal and dual solutions satisfy the following complementary slackness conditions
\begin{equation} \begin{array}{rl}\label{c}
\displaystyle \sum_{i=1}^{m}\pi_i \left(p_i -x\right)=0\quad \text{and}\quad \pi_i\left(r^2 - \|x-p_i\|^2\right) = 0,\quad i=1,\dots,m.
\end{array}\end{equation}

\noindent An important consequence of duality and the complementary slackness conditions, and a well known fact, is the following lemma

\begin{lemma}\label{lem:convhull} 
Consider a ball $B(x,r)$ that covers $\Pp$, and let $\Ss\subseteq\Pp$ be the set of points on the boundary of $B(x,r)$. Then $B(x,r)=MB(\Pp)$ if and only if $x\in\conv(\Ss)$.
\end{lemma}

\noindent As a consequence, the optimal ball is determined by at most $n+1$ affinely independent points of~$\Pp$. Moreover, it is easily proved that the center is the intersection of the bisectors of the facets of $\conv(\Ss)$, the convex hull of $\Ss$.

\medskip

We now define \emph{support set} for the $MB(\Pp)$ problem, an analogous concept to the one of basis for linear programming. 

\begin{definition}[Support set]
	An affinely independent subset $\Ss\subseteq\Pp'\subseteq\Pp$ is a \emph{support set} of $MB(\Pp')$ if $MB(\Ss) = MB(\Pp')$ and $\Ss$ is minimal, in the sense that it does not exist a $\Ss'\varsubsetneq\Ss$ s.t. $MB(\Ss')=MB(\Pp')$.
\end{definition}

\noindent The following lemma is yet another consequence from duality and complementary slackness:

\begin{lemma}\label{lem:suppset}
$\Ss$ is a support set of $\Pp'\subseteq\Pp$ if and only if $x$, the center of $MB(\Pp')$, is in $\ri\conv(\Ss)$, the relative interior of $\conv(\Ss)$. 
\end{lemma}

The dual algorithm developed by Dearing and Zeck \citeyearpar{Dearing09} solves the problem $MB(\Pp)$ by finding a support set of $MB(\Pp)$ the following way: at the beginning of each iteration we have the solution of $MB(\Ss)$, where $\Ss$ is a support set of $MB(\Ss)$; a point $p\in\Pp$ that is not yet covered by $MB(\Ss$) is then selected, and problem $MB(\Ss\cup\{p\})$ is solved by iteratively removing points from the set $\Ss\cup\{p\}$ until a support of $MB(\Ss\cup\{p\})$ is found. The radius of the covering ball strictly increases at each iteration and the algorithm stops when all points have been covered. As a consequence, the algorithm is finite \citep[Theo. 3.2]{Dearing09}.

\medskip

We now review the algorithm in more detail. We denote by $x$ and $r$ the center and radius of the ball at each iteration, and by $\Ss$ the corresponding support set. 

As $\Ss$ and $x$ are updated throughout each iteration, the algorithm maintains the following invariants
\begin{itemize} 
	\item $\Ss$ is affinely independent, 
	\item $x\in\conv(\Ss)$.
\end{itemize}
\noindent As a consequence, dual feasibility and complementary slackness conditions are both always satisfied.

\subsection{Initialization}
The routine starts with a support set $\Ss\subseteq\Pp$ of $MB(\Ss)$, and $x$ and $r$ the solution to $MB(\Ss)$. If such data is not available, the algorithm picks any two points $\{p_{i_1}, p_{i_2}\} \in \Pp$, and solves $MB(\Ss)$ with $\Ss=\{p_{i_1},p_{i_2}\}$:
	$$x = \frac{p_{i_1}+p_{i_2}}2\quad\text{and}\quad r = \|x-p_{i_1}\|.$$
	
\subsection{Iteration}	
	At the beginning of each iteration, we have the minimum covering ball of a support set $\Ss= \{p_{i_1}, . . . , p_{i_s}\}\subseteq \Pp$, whose center and radius are $x$ and $r$ respectively. Moreover $x\in\ri \conv (\Ss)$.
	
	\medskip
	
\begin{description}
	\item[1. \emph{Optimality check}:] If all points of $\Pp$ are covered by $MB(\Ss)$ then $MB(\Ss)=MB(\Pp)$. Otherwise, the algorithm picks a point $\pp\in\Pp$ that is not yet covered. 
	\bigskip
	
	\item[2. \emph{Update $\Ss$}:] If $\Ss\cup\{\pp\}$ is affinely independent then $\Ss=\Ss\cup\{\pp\}$, and the two invariants are maintained for the current $\Ss$ and $x$. 

	If $\Ss\cup\{\pp\}$ is not affinely independent, a point $p_{i_k}\in\Ss$ is dropped from $\Ss$, in such way that $\Ss\setminus\{p_{i_k}\}\cup\{\pp\}$ is affinely independent and $x\in\conv(\Ss\setminus\{p_{i_k}\}\cup\{\pp\})$. Consider $s$ the cardinality of $\Ss$. Point $p_{i_k}$ to leave $\Ss$ is calculated as follows:

	\begin{itemize}
		\item Let $\pi_1,...,\pi_s$ solve  $\displaystyle\sum^s_{j=1} \pi_j p_{i_j} = x\quad\text{and}\quad \sum^s_{j=1} \pi_j=1$\label{upS1};
		\item Let $\omega_1,...,\omega_s$ solve  $\displaystyle\sum^s_{j=1} \omega_j p_{i_j} = -\pp\quad\text{and}\quad \sum^s_{j=1} \omega_j=-1$\label{upS2};		
		\item $p_{i_k}$ is such that $\displaystyle\frac{\pi_k}{-\omega_k} = \min_{j=1,\dots,s} \left\{\frac{\pi_j}{-\omega_j}: \omega_j<0\right\}$\label{upS3}.
	\end{itemize}
	$\Ss$ is now updated: $\Ss=\Ss\setminus \{p_{i_k}\}\cup\{\pp\}$. 
	
	\bigskip
	
	\item[3. \emph{Solution of $MB(\Ss)$}:]  The optimal center of $MB(\Ss)$ can occur in two possible locations:
	
	(a) Either the center is in the interior of $\conv(\Ss)$, in which case it is the intersection of the bisectors of the facets of $\conv(\Ss)$ with $\aff(\Ss)$, the affine space generated by the points in $\Ss$. In this case $\Ss$ is the support set of $MB(\Ss)$.
	
	(b) Or the center is on one of the facets of $\conv(\Ss)$, implying that the support set of $MB(\Ss)$ is a proper subset of $\Ss$.
	
	\medskip
	
	To find the optimal center of $MB(\Ss)$, the algorithm uses a directional search procedure that starts at the current center $x$, and proceeds along a direction $d$. During this search the algorithm either immediately finds the optimal solution to $MB(\Ss)$, or identifies a point from $\Ss$ that is not part of the support set of $MB(\Ss)$, removing it from $\Ss$, and performing a new directional search next. The details of this directional search are described below:
	
	\smallskip
	
	\begin{description}
		\item [i. \emph{The direction $d$}: ] $d$, the direction along which the line search is performed, satisfies the following:
		
		(a) It is on the $(s-1)$-dimensional subspace generated by the points in $\Ss$, $\Sub(\Ss)$, where~$s$ is the cardinality of $\Ss$, that is
		$$u_j^Td = 0,\,\, j=1,...,n-s+1,$$
		\noindent with $\{u_j\}_j$ a basis for $\Null(\Ss)$, the null space of $\Sub(\Ss)$;
		
		(b) It is parallel to the intersection of the bisectors of the facets of the polytope $\conv(\Ss\setminus\{\pp\})$, or, equivalently, it is orthogonal to $\Sub(\Ss\setminus\{\pp\})$, so
		$$(p_{i_j}-p_{i_s})^Td=0,\quad p_{i_j}\in\Ss\setminus\{p_{i_s},\pp\};$$
		
		(c) It ``points towards''  $\,\pp$, in the sense that the distance to $\pp$ from any point on the ray $\ell^+ = \{x+\alpha d: \alpha \geq 0\}$ decreases as $\alpha$ increases:
		$$(p-p_{i_s})^Td=1.$$
		
		\medskip
		
		\item [ii. \emph{Calculating the next iterate}: ] The ray $\ell^+$ intersects both the intersection of the bisectors of the facets of $\conv (\Ss)$, let $\alpha_b$ correspond to that point, and one (or the intersection of several) of those facets, let $\alpha_f$ correspond to that point. Two cases are possible:
		
		\smallskip
		Case 1: $\alpha_b<\alpha_f$, that is, the intersection with the bisectors occurs first. If this is the case, the solution to $MB(\Ss)$ is the point $x+\alpha_b d$. The algorithm goes back to Step 1 with $x=x+\alpha_b d$ and the support set $\Ss$.
	
		\smallskip
		Case 2: $\alpha_b\geq\alpha_f$. In this case the opposite point $p_{i_l}\in \Ss$ to the (or one of the) intersected facet(s) is not part of the support set of $MB(\Ss)$, being therefore removed from $\Ss$: $\Ss = \Ss\setminus\{p_{i_l}\}$. Note that $p_{i_l}$ can never be $p$. The algorithm now returns to the beginning of Step 3 with the new $\Ss$ and $x=x + \alpha_f d$, for a new directional search to solve $MB(\Ss)$.
	\end{description}
\end{description}

\noindent For a full description of the algorithm and its correctness we refer the reader to \citep{Dearing09}.

%% file: 3NewLineSearch.tex
In \citep{Dearing09}, the authors find $\alpha_f$, the point where the ray $\ell^+$ intersects the boundary of $\conv (\Ss)$, by calculating the intersection of the ray with each one of the facets (as many as $n+1$). This is the central reason why, even using efficient updates to the QR factorization of a matrix, their iteration could not have a better complexity than $\bigO(n^3)$, since such updates would need to be done $\bigO(n)$ times. 

\smallskip

We now show how one can find $\alpha_f$ without having to check each facet. The idea consists on projecting $\conv(\Ss)$ and the ray $\ell^+$ orthogonally onto $\aff(\Ss\setminus\{\pp\})$. Recall that $\ell^+$ is perpendicular to $\aff(\Ss\setminus\{\pp\})$, so its projection will be a single point, the projection of $x$. In order to find the intersected facet, we find in which two projected facets of $\conv(\Ss)$ the projection of $x$ fell into. We then calculate the intersection of the ray with the two facets, and the one with smallest $\alpha$ that is non-negative, corresponds to the facet of $\conv(\Ss)$ intersected by the ray.

\medskip
 
\noindent Before we proceed into the details, consider the following notation:
\begin{itemize}
	\item $\Ss=\{p_1, \dots, p_{s-1}, \pp\}$ and $\Ss'=\Ss\setminus\{\pp\}$;
	\item $\C$, the polytope $\conv(\Ss)$, and $\partial \C$ its boundary;
	\item $F_j= \conv(\Ss\setminus\{p_j\})$, the facet of $\C$ opposed to point $p_j\in\Ss$, $j=1,...,s-1$;
	\item $F_0=\conv(\Ss\setminus\{\pp\})$, the facet of $\C$ opposed to~$\pp$;
	\item $\ell=\{x+\alpha d: \alpha\in\RR\}$ and $\ell^+=\{x+\alpha d: \alpha\geq 0\}$;
	\item $\alpha_j$, $j=0,\dots,m$, be the intersection of $\ell$ with facet $F_j$.
\end{itemize}

\medskip

Recall that, at the beginning of Step 3ii, we have:
\begin{itemize}
	\item $\Ss$ is an affinely independent set, and therefore $\pp\not\in\aff(\Ss')$;
	\item $d$ is a direction in $\Sub(\Ss)$, orthogonal to $\aff(\Ss')$, that points towards $\pp$ and passes through the intersection of the bisectors of the facets of $\conv(\Ss)$;
	\item $x$, the current solution, is in $\conv(\Ss)$.
\end{itemize} 


\medskip

Let $x'$ and $p'$ be the orthogonal projections of $x$ and $\pp$ onto $\aff(\Ss')$, respectively, and let
\begin{equation}\label{rep}
x'=\sum_{j=1}^{s-1} \pi_j p_j\text{ s.t. } \sum_{j=1}^{s-1} \pi_j =1\quad\text{and}\quad p'=\sum_{j=1}^{s-1} \omega_j p_j\text{ s.t. } \sum_{j=1}^{s-1} \omega_j=1,
\end{equation}

\noindent be their unique representations as affine combinations of points of $\Ss'$. An important observation is that the projection of all points of $\ell$ onto $\aff(\Ss')$ coincides with $x'$. The two Lemmas that follow will be useful later.

\begin{lemma}\label{lem:1}
 $x'\in\conv(\Ss'\cup\{p'\})$.
\end{lemma}

\noindent Lemma \ref{lem:1} is a direct consequence of the linearity of the projection operator and the fact that $x\in\conv(\Ss)$. 

\begin{lemma}\label{lem:2}
	If $\pi_j<0$ then $\omega_j<0$.
\end{lemma}
\proof{
From Lemma \ref{lem:1} we know that there exists $\beta\geq 0$ and  $\beta'\geq 0$ such that
$$x'=\sum_{j=1}^{s-1}\beta_j p_j+\beta'p'\text{ and } \sum_{j=1}^{s-1}\beta_j+\beta'=1.$$

\noindent Therefore,
$$x' = \sum_{j=1}^{s-1}\beta_jp_j + \beta'\sum_{j=1}^{s-1}\omega_jp_j = \sum_{j=1}^{s-1}(\beta_j+\beta'\omega_j)p_j.$$

\noindent Let $\delta_j = \beta_j+\beta'\omega_j$. Note that $\sum_{j=1}^{s-1} \delta_j= 1$. Since $\Ss$ is an affinely independent set, the representation of $x'$ as an affine combination of the points in $\Ss$ is unique, therefore we must have $\pi_j=\delta_j=\beta_j+\beta'\omega_j$ for all $j$. Thus, if $\pi_j<0$ then $\omega_j<0$. Note that $\pi_j=0$ does not imply $\omega_j\leq 0$.}

\medskip

Consider the general case where the intersection of line $\ell$ with $\partial \C$ is two points $z_1$ and $z_2$. Let $F_{k_1}$ be one of the facets containing $z_1$ and $F_{k_2}$ be one containing $z_2$. Point $z_1$ and/or $z_2$ may be on the intersection of several facets, but knowing one of them suffices. The projection of $z_1$ and $z_2$ onto $\aff(\Ss')$ is $x'$, therefore $x'$ will be written as a unique convex combination of the projections of the points that form $F_{k_1}$ and also as a unique convex combination of the projections of the points that form $F_{k_2}$, and only of those and no other projected facets. In the particular cases when $\ell$ intersects $\partial\C$ on a single point or on an infinite number of points, $x'$ will still be written as a convex combination of the projected points of one of the intersected facets. Lemma \ref{lem:3} proves this fact. Note that the intersection always exists since $x\in\ell\cap\C$. Theorem \ref{theo:main} finds those at most two unique convex representations of $x'$, finding consequently the intersected facets of $\partial \C$ by $\ell$.

\begin{lemma}\label{lem:3}
	Line $\ell$ intersects $F_k$, the facet of $\C$ opposed to $p_k\in\Ss'$, if and only if $x'\in \conv(\Ss'\setminus\{p_k\}\cup\{p'\})$. 
\end{lemma}
\proof{
	We only need to prove that if $x'\in \conv(\Ss'\setminus\{p_k\}\cup\{p'\})$ then $\ell$ intersects $F_k$, since the opposite is trivial as a consequence of the linearity of the projection operator.
	
	Let $B=[p_2-p_1, \dots, p_{s-1}-p_{1}]$. Then $p'=B(B^TB)^{-1}B^T(\pp-p_1)+p_1$. Firstly, we prove that there exists a $ \gamma >0$ s.t. $d = \gamma (\pp-p')$. Clearly, $\pp-p'\in\Sub (\Ss)$, and $\pp-p'$ is orthogonal to $\Sub (\Ss')$ since
	$$B^T(\pp-p') = B^T(\pp-p_1) - B^T(\pp-p_1)=0.$$	
	\noindent Moreover, $d$ and $\pp-p'$ have the same direction since (recall that $(\pp-p_{s-1})^Td>0$)
	$$(\pp-p_{s-1})^T (\pp-p')= \norm{\pp-p'}^2 + (p'-p_{s-1})^T(\pp-p') = \norm{\pp-p'}^2 >0.$$
	\noindent Thus we conclude that there exists a $\gamma>0$ such that $d=\gamma (\pp-p')$.	
 
	 Now, suppose  $x'\in \conv(\Ss'\setminus\{p_k\}\cup\{p'\})$, that is, 
	$$x'=\sum_{j\neq k} \beta_jp_j + \beta'p'\,\,\text{with}\,\,\sum_{j\neq k}\beta_j+\beta'=1\,\,\text{and}\,\, \beta'\geq 0,\,\beta_j\geq 0,\, j=1,...,s-1.$$
	
	\noindent Observe that for any $\alpha$ we have
	\begin{equation*}
		x'+\alpha d = \sum_{j\neq k} \beta_jp_j + \beta'p' + \alpha\gamma(p-p')=\sum_{j\neq k} \beta_jp_j + (\beta'-\alpha\gamma)p' + \alpha\gamma \pp.
	\end{equation*}
		
	\noindent Let $\alpha' = \frac{\beta'}{\gamma}$. We have that $x'+\alpha' d\in\conv(\Ss'\setminus\{p_k\}\cup\{\pp\})$ since $\sum_{j\neq k}\beta_j+\alpha'\gamma=1$.
	This implies that there exists an $\alpha$ such that $x+\alpha d\in F_k$, that is, $\ell$ intersects $F_k$.
}	

\bigskip

Now we show in Theorem \ref{theo:main} how to calculate the intersection of $\ell^+$ and $\partial\C$, that is $\alpha_f$. Recall that
$$\alpha_f = \min_{j=1,\hdots,s-1} \{\alpha_j:\, \alpha_j\geq 0\},$$

\noindent and $\alpha_j$, the intersection of $\ell$ with any $F_j$, can be calculated by finding a basis $\{w_1,\dots,w_{n-(s-1)+1}\}$ for the null space of $\Sub(\Ss'\setminus\{p_k\}\cup\{\pp\})$, and
\begin{equation}\label{alpha}
\alpha_j = \frac{(\pp-x)^Tw_i}{d^Tw_i},\text{ for any }\, i=1,...,n-s+2, \text{ s.t. } d^Tw_i\neq 0.
\end{equation}
	
\smallskip 

\begin{theorem}\label{theo:main}
	
	Consider the representations (\ref{rep}) of $x'$ and $p'$. Let $F_k$, the facet opposed to $p_k\in\Ss'$, be the one first intersected by~$\ell^+$, and let $\alpha_k$ be the value of $\alpha$ at which the intersection occurs. To find $p_k$ and $\alpha_k$, there are two possible cases:
	\begin{itemize}
		\item {Case 1.} If there is a $k=1,...,m$ such that $\pi_k=\omega_k=0$, then $p_k$ is the point opposed to the facet intersected first, and $\alpha_f=\alpha_k=0$.
		\smallskip
		\item{Case 2.} Suppose case 1 does not hold. First, find ${k_1}$ such that
		\begin{equation}\label{pk1}
		\frac{\pi_{k_1}}{\omega_{k_1}} = \min_{j=1,...,s-1} \left\{\frac{\pi_j}{\omega_j}: \pi_j\geq 0, \omega_j>0\right\},
		\end{equation}
		
		\noindent then, find $\alpha_{k_1}$ as in (\ref{alpha}). Let $\J:=\{j: \pi_j\leq 0, \omega_j<0\}$. If $\J\neq\emptyset$,  find ${k_2}$ such that
		\begin{equation}\label{pk2}
		\frac{\pi_{k_2}}{\omega_{k_2}} = 
		\max_{j=1,...,s-1} \left\{\frac{\pi_j}{\omega_j}: \pi_j\leq 0, \omega_j<0\right\} ,
		\end{equation}
		\noindent and find $\alpha_{k_2}$ as in (\ref{alpha}). If $\J=\emptyset$ simply consider $\alpha_{k_2}=-\infty$.
		
		\noindent The facet first intersected by $\ell^+$, $F_k$, is such that
		$$k=\arg\min_{j=k_1, k_2}\{\alpha_j:\,\alpha_j\geq 0\},$$
		\noindent  and $\alpha_f=\alpha_k$.
	\end{itemize}
\end{theorem}

\proof{
{\emph{Case 1.}} Suppose there is a $k\in\{1,\dots,m\}$ such that $\omega_k=\pi_k=0$. Then, since $x=x'+\delta d$, for some $\delta\geq 0$, and the fact that there is a $\gamma>0$ such that $d=\gamma (p-p')$ (see the proof of Lemma \ref{lem:3}) we have:
\begin{eqnarray*}\label{xx}
	x &=&  x'+\delta\gamma (p-p')=\sum_{\substack{j=1\\ j\neq k}}^{s-1} \pi_jp_j + \delta \gamma \left(p- \sum_{\substack{j=1\\ j\neq k}}^{s-1}  \omega_jp_j \right)=\sum_{\substack{j=1\\ j\neq k}}^{s-1} \left(\pi_j-\delta\gamma\omega_j\right)p_j +\delta \gamma p,
\end{eqnarray*}
			
\noindent and this representation of $x$ as a convex combination of $\Ss$ is unique, since $\Ss$ is affinely independent. Consequently, $\pi_j-\delta\gamma\omega_j\geq 0$ and $\delta\gamma\geq 0$, concluding that $x\in\conv(\Ss'\setminus\{p_k\}\cup\{p\})\equiv F_k$. Therefore $\ell^+$ intersects $F_k$ at $\alpha_f=\alpha_k=0$.

		
\smallskip 

\noindent \emph{Remark: The case $\pi_k=\omega_k=0$ needed to be treated separately, since in such case~$F_k$ is perpendicular to $F_0$, and so formula (\ref{alpha}) could not be applied (we would have~$0/0$). }
		
\smallskip
		
{\emph{Case 2.}} Since $x'\in\conv(\Ss'\cup\{p'\})$, from the proof of Lemma~\ref{lem:2}, there exists $\beta'\geq 0$ such that

\begin{equation}\label{convxp}
	x'=\sum_{j=1}^{s-1}\beta_j p_j+\beta' p'= \sum_{j=1}^{s-1}(\pi_j-\beta'\omega_j)p_j+\beta' p'.
\end{equation}

\noindent Formula (\ref{convxp}) gives all possible ways to represent $x'$ as a convex combination of $\Ss'\cup\{p'\}$ as a function of $\beta'$. We are now interested in knowing the minimum and maximum values of $\beta'$, $\beta'_{min}$ and $\beta'_{max}$ respectively.

Suppose $\pi_j\geq 0$ for all $j=1,...,s-1$. We have that $x'\in\conv(\Ss'\cup\{p'\})$ and so $\beta_{min} = 0$. It is easy to see that $\beta_{max}=\frac{\pi_{k_1}}{\omega_{k_1}}$ as in (\ref{pk1}), and that any $\beta\in[0, \beta_{max}]$ yields $\pi_j-\beta\omega_j\geq 0$. We then conclude that, when all $\pi_j\geq 0$ we have~$x'\in \conv(\Ss'\setminus\{p_{k_1}\}\cup\{p'\})$, and so from Lemma \ref{lem:3} we have that $\ell$ intersects~$F_{k_1}$. Note that $\beta_{max}$ may be $0$. When $\beta_{max}>0$, observe that there is no other $\beta\in]0, \beta_{max}[$ such that $\pi_j-\beta \omega_j= 0$, so there is no other way to write $x'$ as a convex combination of $p'$ and $s-2$ points of $\Ss'$.

\smallskip

Now consider that there exists $\pi_j<0$. Then $x'\not\in \conv(\Ss')$ and therefore $\beta_{min}>0$. Since $\beta'$ must be positive, in order to have a convex combination in (\ref{convxp}), $\beta'$ must satisfy the following conditions:
\begin{equation}\label{neg}
\beta'\geq\frac{\pi_j}{\omega_j},\quad\forall j: \, \omega_j<0,\, \pi_j\leq 0\,\Rightarrow\,\beta'\geq\frac{\pi_{k_2}}{\omega_{k_2}},
\end{equation}
\begin{equation}\label{pos}
\beta'\leq\frac{\pi_j}{\omega_j},\quad\forall j: \, \omega_j>0,\, \pi_j\geq 0\,\Rightarrow\,\beta'\leq\frac{\pi_{k_1}}{\omega_{k_1}}.
\end{equation}

\noindent For $k_1$ and $k_2$ as in (\ref{pk1}) and (\ref{pk2}) respectively. The conditions above are feasible since there must be such a $\beta'$, and because of Lemma \ref{lem:2}. This allow us to conclude that $\frac{\pi_{k_2}}{\omega_{k_2}}\leq\frac{\pi_{k_1}}{\omega_{k_1}}$ and $\beta_{min} = \frac{\pi_{k_2}}{\omega_{k_2}}$ and $\beta_{max} = \frac{\pi_{k_1}}{\omega_{k_1}}$. Finally we conclude that $x'\in \conv(\Ss'\setminus\{p_{k_1}\}\cup\{p'\})$ and $x'\in \conv(\Ss'\setminus\{p_{k_2}\}\cup\{p'\})$, with $k_1\neq k_2$. So $\ell$ intersects both facets $F_{k_1}$ and $F_{k_2}$.
}


%% file: 4Implementation.tex
Recall the QR factorization of some matrix $A$ of size $n\times m$ into the product of matrices $Q$ and $R$, where $Q_{n\times n}$ is an orthogonal matrix and $R_{n\times m}$ is upper triangular. For our purposes we will consider $n\geq m$. If $A$ has full column rank, then the diagonal of $R$ is non-zero, and $Q$ can be partitioned in $[V_{n\times m}\, \,U_{n\times (n-m)}]$ such that the columns of $V$ form a basis of $\Span(A)$, the column range of $A$, and the columns of $U$ form a basis to $\Null(A)$, the null space of $A$. Different algorithms are available to find a QR factorization of a matrix, but in terms of computational work, and for a general matrix $n\times n$, they all need $\bigO(n^3)$ steps \citep[\S 5.2]{Golub96}. However, the QR factorization of $A$ can be ``recycled'' and used to calculate the QR factorization of a matrix obtained from $A$ by either rank-one changes, appending a row or column to $A$, or deleting a row or column from $A$ \citep[\S 12.5]{Golub96}. This is accomplished by using \emph{Givens rotations}, and, in the case when $m=n$, such procedure needs $\bigO(n^2)$ steps.

\medskip

We now describe in detail how one can implement the algorithm taking advantage of  the QR factorization updates. At the beginning of each iteration, we have the QR factorization of the $n\times s$ matrix $S$, whose columns are the points in $\Ss$, the support set found in the previous iteration:
$$S =\left[\begin{array}{cccc} p_{i_1} &p_{i_2}&\dots&  p_{i_s}\end{array}\right].$$
Let $Q_S$ and $R_S$ be the matrices of the QR factorization of $S$, which were inherited from the previous iteration. The first iteration is the only time a QR factorization from scratch is performed. 

\subsection{Update $\Ss$ procedure}\label{subsec:app1}

The first step of this phase is to check whether $\Ss\cup \{\pp\}$ is affinely independent. If it is not, the next step is to find $p_{i_k}\in \Ss$, such that $\Ss\setminus\{p_{i_k}\} \cup \{\pp\}$ is affinely independent.

\smallskip

Consider matrices $B$ and $\bar{B}$ as follows:
\begin{equation}\label{mat:B}
B = \left[\begin{array}{ccc}
p_{i_1}& \dots & p_{i_s}\\
1 & \dots &1
\end{array}\right], \quad\quad \bar{B} = \left[\begin{array}{cccc}
p_{i_1}& \dots & p_{i_s} & \pp\\
1 & \dots &1 & 1
\end{array}\right].
\end{equation}

\noindent $B$  is full column rank since $\Ss$ is affinely independent. When $\Ss$ has $n+1$ points we automatically know that $\Ss\cup \{\pp\}$ is affinely dependent. Consider then $s\leq n$, which implies that $\bar{B}$ has at least as many columns as rows. $\Ss\cup \{\pp\}$ is affinely independent if and only if $\bar{B}$ is full column rank. This can be checked by looking at the element in position $(s+1, s+1)$ of $R_{\bar{B}}$, from the QR factorization of $\bar{B}=Q_{\bar{B}}R_{\bar{B}}$: if it is zero or close to zero (to account for precision errors) then $\Ss\cup \{\pp\}$ is affinely dependent, otherwise it is affinely independent.


Matrix $B$ can be obtained by inserting a row of ones in $S$, and $\bar{B}$ by then inserting the column $\binom{\pp}{1}$ in $B$, thus the corresponding QR factorizations can be efficiently computed from the QR factorization of $S$ and of $B$ respectively.

\smallskip

If $\Ss\cup \{\pp\}$ is affinely independent, we leave matrix $S$ as is. Otherwise, the linear systems are solved
\begin{equation}\label{ls:B}
B\pi = \binom{x}{1},\quad\text{and}\quad B\omega = \binom{-\pp}{-1},
\end{equation}
\noindent by reducing them to linear systems with upper triangular matrices, using the QR factorization of $B$:
\begin{equation}\label{ls:RB}
R_B\pi = Q_B^T\binom{x}{1},\quad\quad\quad R_B\omega = Q_B^T\binom{-\pp}{-1},
\end{equation}
\noindent which can be solved using \emph{Back Substitution} and performed in $\bigO(n^2)$ \citep[\S 3.1]{Golub96}. Now $\Ss$ is updated, $\Ss=\Ss \setminus \{p_{i_k}\}\cup\{\pp\}$, for some $p_{i_k}$, and at this point we get the new $Q_S$ and $R_S$, the QR factorization of the new matrix $S$, obtained from the previous one by removing the $k$-th column. We do not add $\pp$ yet to matrix $S$, we leave that to the very end of the iteration, when $\pp$ is finally covered. 
%
%

\subsection{Solving MB($\Ss$): calculating $d$}\label{subsec:app2}

Consider that at the beginning of this phase $S= \left[p_{i_1}\, \dots \, p_{i_{s}}\right]$ and let $C$ be the following matrix:		
\begin{equation}	\label{mat:C}C = \left[\begin{array}{cccc}
p_{i_1}-p_{i_s}& \dots & p_{i_{s-1}}-p_{i_s} & \pp - p_{i_s}
\end{array}\right].
\end{equation}

\noindent Let $C = Q_{C}R_{{C}}$ be the QR factorization of ${C}$, obtained by updating the QR factorization of $S=Q_SR_S$ twice, since $C$ can be obtained from $S$ by adding two rank one matrices
$${C} = S + (\pp-p_{i_s})e_n^T - p_{i_s}1_n^T.$$
\noindent with $e_j\in\RR^n$ is the vector with $1$ in the $j$-th entry and all the other zero. ${C}$ is full column rank since $\Ss\cup\{\pp\}$ is affinely independent, so we can partition $Q_{C}=[V\,\,  U]$, where $V$ is a $n\times s$ matrix whose columns are a basis to $\Span(C)$, and $U$ is a $n\times(n-s)$ matrix whose columns are a basis to $\Null(C)$. Finally, to find $d$ we need to solve the linear system

\begin{equation*}\label{ls:C}
\left[\begin{array}{c}
C^T\\
U^T
\end{array}\right] d=e_{s}\quad\Longleftrightarrow\quad
\left[\begin{array}{ccc}
& R_C^T&\\
\hline
0_{(n-s)\times s} & \vline &I_{n-s}
\end{array}\right]Q_C^T d = e_s,\end{equation*}

\noindent where the latter is a lower triangular system.
	
		
\subsection{Solving MB($\Ss$): calculating the next iterate}\label{subsec:app3}
		
The first step now is to calculate $x'$ and~$\pp'$, the orthogonal projection of $x$ and $\pp$ respectively onto~$\aff(\Ss')$, where $\Ss'=\Ss\setminus\{\pp\}$. Let $D$ be the following matrix:
	
\begin{equation*}\label{mat:D}
D = \left[\begin{array}{ccc}
p_{i_1}-p_{i_s}& \dots & p_{i_{s-1}}-p_{i_s}
\end{array}\right].
\end{equation*}

\noindent And let $Q_D$ and $R_D$ be the matrices of the QR factorization of $D$, which can be obtained easily from the QR factorization of matrix $C$ (\ref{mat:C}), since $D$ is obtained from $C$ by removing the last column. Let $V$ be the matrix with the first $s-1$ columns of $Q_D$, which form an orthogonal basis to $\Sub(\Ss')$. Consider $\aff(\Ss')=p_{i_1}+\Sub(\Ss')$. The projections $x'$ and~$\pp'$ can be calculated the following way:

$$x' = VV^T(x-p_{i_1})+p_{i_1},\quad \quad \pp' = VV^T(\pp-p_{i_1})+p_{i_1}.$$

\noindent Now, two linear systems need to be solved
\begin{equation}\label{ls:B2}
B\pi = \binom{x'}{1}\quad\text{and}\quad B\omega = \binom{\pp'}{1},
\end{equation}

\noindent where $B$ is the same matrix as in (\ref{mat:B}) but with the points of $\Ss'$, so (\ref{ls:B2}) are solved following the same steps used to solve (\ref{ls:B}).

After finding $k_1$ and $k_2$, as in (\ref{pk1}) and (\ref{pk2}) respectively, we need to calculate $\alpha_{k_1}$ and $\alpha_{k_2}$. To calculate~$\alpha_{k_1}$, we need a basis for the null space of $\Sub(\Ss'\setminus\{p_{k_1}\}\cup \{\pp\})$, in order to apply formula (\ref{alpha}). If $k_1<s$, $\Sub(\Ss'\setminus\{p_{i_k}\}\cup \{\pp\})\equiv \Span(F)$ with
$$F=\left[\begin{array}{ccccccc}\label{mat:F}
p_{i_1}-p_{i_s}& \dots &p_{i_{k_1-1}}-p_{i_s}&p_{i_{k_1+1}}-p_{i_s}&\dots& p_{i_{s-1}}-p_{i_s} & \pp - p_{i_s}
\end{array}\right].$$

\noindent Matrix $F$ can be obtained from $C$ (\ref{mat:C}) by deleting its $k_1$-th column. On the other hand, if $k_1=s$ then	$\Sub(\Ss'\setminus\{p_{i_k}\}\cup \{\pp\})\equiv \Span(F)$ with
$$F=\left[\begin{array}{ccc}\label{mat:F2}
p_{i_1}-\pp& \dots & p_{i_{s-1}}-\pp
\end{array}\right],$$
\noindent and $F$ can be obtained from $C$ by deleting the last column and then adding the rank one matrix $(p_{i_{s}}-\pp)1_n^T$. Therefore, in both cases, the QR factorization of $F$ can be obtained by updating the QR factorization of~$C$. A basis for the null space of $\Sub(\Ss'\setminus\{p_{k_1}\}\cup \{\pp\})$ is then formed by the  last $n-s+1$ columns of~$Q_F$. The value $\alpha_{k_2}$ is calculated in an analogous way.


\medskip

The only time a $QR$ factorization of a matrix is calculated from scratch is at the beginning of the algorithm. Then, at each iteration, a constant number of ``QR updates'' are performed to matrices with $n$ rows and at most $n$ columns, so such updates take $\bigO(n^2)$ steps. This results in an iteration that is done in quadratic time.

%% file: 5Results.tex
In order to understand in practice the effect of the new directional procedure together with implementation described in section \ref{sec:Implementation}, we implemented the algorithm in MATLAB. In order to compare it to the original version from Dearing and Zeck \citeyearpar{Dearing09}, we re-wrote their algorithm also using QR updates and implemented it too. Our experiments were conducted using MATLAB R2014a (version 8.3) on a PC with an Intel Core i5 2.30 GHz processor, with 4 GB RAM. Tables \ref{tab:1} and \ref{tab:2} show the average running times of the two versions of the algorithm on instances with $1000$ and $10000$ points, respectively, drawn uniformly at random from the unit cube. 
\begin{table}
\caption{Average time in seconds for datasets with $m=1000$ points in variable dimension $n$ uniformly sampled in a unit cube.}
\begin{center}\begin{tabular}{rrrrr}
\hline\noalign{\smallskip}
\multicolumn{2}{c}{Problem} & &\multicolumn{2}{c}{Time in seconds}\\
\noalign{\smallskip}\cline{1-2}\cline{4-5}\noalign{\smallskip}
$n$ &  $m$ && D\&Z original & D\&Z new\\
\noalign{\smallskip}\hline\noalign{\smallskip}
50 & 1000 && 0.15 &	0.32		 \\
100 & 1000 && 0.42 & 0.96 \\
500 & 1000 && 32.52 &	20.13 \\	
1000 & 1000 && 140.45 & 57.50	\\
5000 & 1000 && 8852.20 & 887.15 \\
\noalign{\smallskip}\hline
\end{tabular}\end{center}\label{tab:1}  
\end{table}

\begin{table} 
\caption{Average time in seconds for datasets with $m=10000$ points in variable dimension $n$ uniformly sampled in a unit cube.}
\begin{center}\begin{tabular}{rrrrr}
		\hline\noalign{\smallskip}
		\multicolumn{2}{c}{Problem} & &\multicolumn{2}{c}{Time in seconds}\\
		\noalign{\smallskip}\cline{1-2}\cline{4-5}\noalign{\smallskip}
		$n$ &  $m$ && D\&Z original & D\&Z new\\
		\noalign{\smallskip}\hline\noalign{\smallskip}
		50 & 10000 && 2.48 &	2.95		 \\
		100 & 10000 && 5.03 & 5.97 \\
		500 & 10000 && 70.04 &	42.55 \\	
		1000 & 10000 && 267.68 & 114.88	\\
		5000 & 10000 && 17044.10 & 1463.20 \\
		\noalign{\smallskip}\hline
	\end{tabular}\end{center}\label{tab:2} 
\end{table}
	
We observed that for smaller dimensions the original algorithm is slightly faster, which is expected. It is with larger dimensions that we observe that the new version of the algorithm with the changes we proposed in section \ref{sec:LineSearch} is considerably faster. This change naturally does not affect the number of iterations, but only the computational work of each iteration.